\documentclass[11pt,letterpaper]{article}
\usepackage[margin=1in]{geometry}
\usepackage[T1]{fontenc}
\usepackage{xcolor}
\usepackage{lmodern,microtype}
\usepackage{amsmath,amsthm,amssymb}
\usepackage{enumerate}
\usepackage{cite}
\usepackage{tikz}
\usetikzlibrary{calc,positioning,decorations.pathreplacing}
\usepackage[colorlinks,allcolors=blue,bookmarksdepth=2]{hyperref}
\usepackage[capitalise,noabbrev]{cleveref}
\numberwithin{equation}{section}
\newtheorem{theorem}{Theorem}[section]
\newtheorem{proposition}[theorem]{Proposition}
\newtheorem{lemma}[theorem]{Lemma}
\theoremstyle{definition}

\newtheorem{remark}[theorem]{Remark}
\newcommand{\softO}{\widetilde O}
\newcommand{\ZZ}{\mathbb Z}
\newcommand{\CC}{\mathbb C}
\newcommand{\iu}{\mathrm{i}}
\newcommand{\tool}[1]{\cref{sec:known}\eqref{#1}}

\title{Finding Integer Roots of Polynomials in Nearly Optimal Time}
\date{}
\author{%
  Yansong Feng\thanks{%
    Academy of Mathematics and Systems Science,
    Chinese Academy of Sciences.
    Email: \href{mailto:fengyansong@amss.ac.cn}
                {\texttt{fengyansong@amss.ac.cn}}.%
  }
  \and
  Itamar Nir\thanks{%
    Tel Aviv University.
    Email: \href{mailto:nirgottlieb@mail.tau.ac.il}
                {\texttt{nirgottlieb@mail.tau.ac.il}}.
    This research was funded by the European Union (ERC, EACTP, 101142020). Views and opinions expressed are however those of the author(s) only and do not necessarily reflect those of the European Union or the European Research Council Executive Agency. Neither the European Union nor the granting authority can be held responsible for them.%
  }
}

\begin{document}
\maketitle

\begin{abstract}
Let $f\in\mathbb Z[x]$ be a polynomial of degree $n$ with integer coefficients
of absolute value less than $2^b$. Using Pan's approximate factorization
theorem, we give a deterministic algorithm that finds all distinct
integer roots of $f$ in $\widetilde O(nb)$ bit operations. This matches the
bit size of the input up to logarithmic factors and does not require
$f$ to be square-free. Since we were unable to verify parts of the
running time analysis of Pan's theorem, we also give a self-contained
deterministic algorithm with running time
$\widetilde O(nb+\min(n,b)^3b^2)$, which is $\widetilde O(nb)$ whenever $b\leq
n^{1/4}$. This answers a question of Harvey and Hittmeir (2022), who
asked whether integer roots can be found deterministically in
quasilinear time, conditionally on Pan's theorem, and unconditionally
for $b\leq n^{1/4}$. The fastest previous deterministic algorithm, due
to Nir, takes $\widetilde O(n^{3/2}b)$ bit operations and works only for
square-free polynomials.

Instead of $p$-adic lifting, we combine approximate complex root finding
with exact polynomial evaluation. We first approximate the complex roots
of the polynomial formed by the top $\min(n,b)+1$ coefficients of $f$.
Every integer root $r$ of $f$ with $|r|\geq 2$ lies within distance $2$
of some complex root of this truncated polynomial. Rounding the real
parts of the approximate roots therefore gives a list of $O(\min(n,b))$
candidate integers that is guaranteed to contain all integer roots of
$f$. Although a single candidate can have up to $b$ bits, a bound via
the Mahler measure shows that the bit lengths of all candidates sum to
$O(b)$. This allows us to test all candidates by exact multipoint
evaluation in $\widetilde O(nb)$ bit operations.
\end{abstract}

{\small\tableofcontents}

\clearpage

\section{Introduction}
\label{sec:introduction}

Let $n,b\geq1$ be integers, and let
\[
 f(x)=a_nx^n+\cdots+a_1x+a_0\in\ZZ[x]
\]
be a polynomial of degree $n$ with $a_n\neq0$ and $|a_i|<2^b$ for all
$i$. The goal is to compute the set of integer roots $\{r\in\ZZ : f(r)=0\}$.
The coefficients of $f$ are given in binary, and the distinct integer
roots are to be output in binary without multiplicities. Running time is measured in bit operations, and we use the fact that
$k$-bit integers can be multiplied and divided with remainder in
$\softO(k)$ bit operations~\cite{HvdH21}. Throughout the paper, $\softO$ hides factors
polynomial in $\log n$ and $\log b$.

A dense polynomial of degree $n$ with $b$-bit coefficients has bit size
$\Theta(nb)$, so any algorithm that reads the whole input must take at
least $\Omega(nb)$ bit operations in the worst case. Our main result
shows that, using Pan's approximate factorization theorem, this bound can
be matched up to logarithmic factors.

\begin{theorem}\label{thm:integer-roots}
Given a polynomial $f\in\ZZ[x]$ of degree $n\geq 1$ with coefficients of
absolute value less than $2^b$, one can deterministically compute all
distinct integer roots of $f$
\begin{enumerate}[\rm(a)]
 \item in $\softO(nb+\min(n,b)^3b^2)$ bit operations, and
 \item in $\softO(nb)$ bit operations, using Pan's approximate
 factorization theorem (Theorem~\ref{thm:pan}).
\end{enumerate}
\end{theorem}

Part~(b) relies on the running time bound of Pan's
algorithm~\cite{Pan02}, as stated by Mehlhorn, Sagraloff, and
Wang~\cite[Theorem~2]{MSW15}. While preparing this paper we were unable
to verify parts of Pan's running time analysis; see
Remark~\ref{rem:pan}. Part~(a) does not use Pan's theorem, and its proof
is self-contained: besides Appendix~\ref{sec:elementary}, it uses only
Lemmas~\ref{lem:landau}, \ref{lem:leading-part},
and~\ref{lem:multipoint}.

Theorem~\ref{thm:integer-roots}(b) answers a question of Harvey and
Hittmeir~\cite[Remark~A.8]{HH22}, who asked whether the integer roots of
a square-free polynomial can be found deterministically in quasilinear
time, that is, in $O(n^{1+\epsilon}b^{1+\epsilon})$ bit operations.
Part~(a) answers it unconditionally when $b\leq n^{1/4}$, since then
$\min(n,b)^3b^2=b^5\leq nb$. In fact, neither algorithm needs $f$ to be
square-free. We discuss their question, and the obstructions they pointed out, in \cref{sec:related}. The main construction in this paper was generated with the assistance of AI models and verified by the authors (see~\hyperref[sec:ai-methodology]{AI methodology} for details).

\subsection{Overview of the algorithm}
\label{sec:overview}

Previous deterministic algorithms for this problem are based on $p$-adic
lifting (see \cref{sec:related}). Here we take a complex-analytic route
instead. Using fast algorithms
for approximate polynomial factorization~\cite{MSW15,Pan02}, we compute
rough complex approximations to the roots of a truncated polynomial,
take the integers close to their real parts, and check the resulting
candidates by evaluating $f$ exactly. In detail, the algorithm runs in
four steps:

\begin{enumerate}[\bf 1.]
\item Let $d=\min(n,b)$ and form the truncated polynomial
 $g(x)=a_nx^d+a_{n-1}x^{d-1}+\cdots+a_{n-d}$.
\item Compute complex numbers $\beta_1,\ldots,\beta_d\in\CC$ such that
 every complex root of $g$ has distance less than $1/8$ to at least one
 $\beta_i$, and
 \[
  \sum_{i=1}^d\log_2\max(1,|\beta_i|)
  \leq b+\tfrac12\log_2(d+1)+2.
 \]
\item For each $i$, pick an integer $z_i$ nearest to
 $\operatorname{Re}\beta_i$, and form the candidate set
 $\mathcal C=\{0,\pm1\}\cup\{z_i+j : 1\leq i\leq d,\ |j|\leq3\}$.
\item Evaluate $f(c)$ exactly for every $c\in\mathcal C$ and output those
 $c$ for which $f(c)=0$.
\end{enumerate}

Since Step~4 evaluates $f(c)$ exactly, any integer returned by the
algorithm is certainly a root of $f$. We need to show that $\mathcal C$
contains all integer roots of $f$, and that Steps~2 and~4 can be
implemented in $\softO(nb)$ bit operations. Before giving the proofs, let
us briefly explain the two main difficulties and how we handle them.

The first issue is that $n$ can be much larger than $b$. Approximating
the roots of a degree-$n$ polynomial to the required accuracy using
Pan's algorithm~\cite{Pan02} needs a working precision of roughly
$n\log n+b$ bits and takes $\softO(n(n+b))$ bit operations, which is too
slow when $b\ll n$. Step~1 gets around this by dropping the low-order
terms of $f$ whenever $b<n$. As we show in Lemma~\ref{lem:leading-part},
any integer root $r$ of $f$ with $|r|\geq2$ lies within distance $2$ of
some complex root of the top part $g$. The reason is simple: when $b<n$
we have $d=b$, and writing $f(x)=x^{n-d}g(x)+h(x)$ with $\deg h<n-d$
gives $|g(r)|=|h(r)|/|r|^{n-d}<2^b=2^d$ at any integer root $r$ with
$|r|\geq2$. On the other hand, if $r$ had distance at least $2$ from
every root $\alpha_j$ of $g$, then
$|g(r)|=|a_n|\prod_{j=1}^d|r-\alpha_j|$ would be at least $2^d$. This
allows us to work in degree $d=\min(n,b)$, bringing the cost of Step~2
down to $\softO(d(b+d))\subseteq\softO(nb)$. With root approximation
error $1/8$ in Step~2, the distance from an integer root $r$ ($|r|\geq
2$) to the nearest integer $z_i$ is less than $2+\tfrac18+\tfrac12<3$,
so checking shifts $|j|\leq3$ in Step~3 is enough to catch $r$.


The second issue is the cost of testing the candidates. The set $\mathcal
C$ has at most $7d+3$ elements, and a single candidate can be as large
as $\Theta(2^b)$. If all $O(d)$ candidates had $\Theta(b)$ bits, their
total bit length would be $\Theta(db)$, and evaluating $f$ at all of
them would take roughly $n\cdot db$ bit operations (which is $n^2b$ when
$b\geq n$). However, the roots of $g$ cannot all be large at the same
time. Using Landau's inequality for the Mahler measure together with
Lemma~\ref{lem:perturbation}, we show that the approximations $\beta_i$
satisfy
\[
 \sum_{i=1}^d\log_2\max(1,|\beta_i|)\leq b+O(\log d).
\]
This implies that the total bit length of the candidates,
\[
 S=\sum_{c\in\mathcal C}\bigl\lceil\log_2(2+|c|)\bigr\rceil,
\]
is $O(b)$. Lemma~\ref{lem:multipoint} then evaluates $f$ on $\mathcal C$
in $\softO(nb+nS)=\softO(nb)$ bit operations: it groups the candidates by
bit length, splits $f$ into blocks, evaluates the blocks by multipoint
evaluation, and combines the block values at each
candidate. For related work on fast approximate multipoint evaluation,
see~\cite{KS16,Mor22}.

\paragraph{Without Pan's theorem.} For part~(a) of
Theorem~\ref{thm:integer-roots}, we replace Steps~2 and~3 by an
elementary bisection search on the real line
(Appendix~\ref{sec:elementary}). It discards an interval as soon as a
simple test, based on Pellet's theorem and Graeffe's root squaring as
in~\cite{BSSY18}, certifies that no root of $g$ lies within distance~$2$
of the interval. The search returns $O(d)$ candidates whose bit lengths
again sum to $O(b)$, and it takes $\softO(d^3b^2)$ bit operations.

\subsection{Related work}
\label{sec:related}

\paragraph{The $p$-adic approach.}
The classical way to find integer and rational roots of polynomials over
$\ZZ$ is $p$-adic lifting. Loos~\cite{Loo83} analyzed a Newton--Hensel
algorithm for the rational roots of a square-free integer polynomial,
with running time $O(n^3L^2)$ in classical arithmetic, where $L$ is the
length of $\|f\|_1$.
A root-finding algorithm along the same lines appears in the book of von
zur Gathen and Gerhard~\cite[Theorem~15.21]{vzGG13}. Harvey and
Hittmeir~\cite{HH22} analyzed a deterministic version of this method, as
a subroutine of their deterministic algorithm for finding $r$-power
divisors of integers. They obtained a bound of
$O(n^{2+\epsilon}b^{1+\epsilon})$ bit operations for $b\geq n$ and
general $f$~\cite[Theorem~1.3]{HH22}, as well as the same bound for
square-free $f$~\cite[Proposition~A.6]{HH22}. They asked whether the
latter can be improved to quasilinear~\cite[Remark~A.8]{HH22}. If
randomization is allowed, the $p$-adic method finds the integer roots of
a square-free $f$ in $\softO(nb)$ bit operations; this is sketched
in~\cite[Remark~A.8]{HH22} and proved in~\cite{Nir26}.
Recently, Nir~\cite{Nir26} improved the deterministic
complexity to $\softO(n^{3/2}b)$ for square-free $f$, which appears to
be the best deterministic bound known prior to our work. Nir's algorithm
works with a small prime $p$ for which $f\bmod p$ is not necessarily
square-free, but has repeated roots of small total multiplicity. It
computes the first $b$ base-$p$ digits of the $p$-adic roots of $f$
using a weighted lifting tree. These digit strings are candidates for
being integer roots of $f$. It then checks the resulting at most $n$
candidates with a deterministic verification procedure taking
$\softO(nb+\min(n^2,nb^2))$ bit operations.

\paragraph{Obstructions for the $p$-adic approach.}
Harvey and Hittmeir pointed out two obstructions to a deterministic
quasilinear algorithm along these lines~\cite[Remark~A.8]{HH22}, and a
third one arises when $f$ is not square-free. None of them arises in the
complex-analytic approach.
\begin{enumerate}[\bf 1.]
\item \emph{Choosing the prime.} The $p$-adic method needs a prime $p$
 for which $f\bmod p$ is square-free, so that every root modulo $p$ has a
 unique Hensel lift. Almost all primes have this property, but it is
 hard to predict which ones, and the deterministic algorithm
 of~\cite{HH22} tests every prime up to $O(nb)$.
\item \emph{Verifying the candidates.} Hensel lifting produces up to $n$
 candidates, each a residue modulo $p^k$ with $\Theta(b)$ bits, and one
 still has to check which of them are actual roots of $f$. Computing the
 values $f(c)$ directly is too slow: each of them can have $\Theta(nb)$
 bits, so just writing all of them down takes $\Theta(n^2b)$ bit
 operations. Nir~\cite{Nir26} verifies the
 candidates deterministically in $\softO(nb+\min(n^2,nb^2))$ bit
 operations, which is still not quasilinear. In our algorithm, the
 candidates are obtained by rounding our approximations of the complex
 roots of the truncated polynomial $g$, and we prove that their bit
 lengths sum to $O(b)$. This makes the verification much simpler:
 Lemma~\ref{lem:multipoint} evaluates $f$ at all candidates in
 $\softO(nb)$ bit operations (see the proof of
 Theorem~\ref{thm:integer-roots}).
\item \emph{Computing the square-free part.} For general $f$, Harvey and
 Hittmeir first replace $f$ by $f/\gcd(f,f')$, and no deterministic
 quasilinear algorithm is known for this step~\cite[Remark~A.10]{HH22}.
 The algorithm of~\cite{Nir26} also assumes that $f$ is square-free. Our
 algorithm makes no assumption on the multiplicities of the roots of $f$
 and never computes $\gcd(f,f')$.
\end{enumerate}

\paragraph{Complex root finding.}
Integer roots can also be read off from sufficiently accurate
approximations of the real or complex roots. Loos~\cite{Loo83} compared
his $p$-adic algorithm with an approach of Collins based on real root
isolation, which takes $O(n^4L^2)$ in classical arithmetic. Isolating
all complex roots of an integer polynomial takes
$\softO(n^3+n^2b)$ bit operations with the algorithm of Becker,
Sagraloff, Sharma, and Yap~\cite{BSSY18}, which is comparable to the
bound of Pan~\cite{Pan02}; see also~\cite{MSW15}. Moroz~\cite{Mor22}
isolates the roots of a square-free polynomial in
$\softO(n(b+\log\kappa))$ bit operations, where $\kappa$ is the
condition number of its roots. For square-free integer polynomials,
Mahler's bound gives $\log\kappa=O(nb+n\log n)$~\cite[full version, Proposition~5]{Mor22},
and hence a worst-case bound of $\softO(n^2b)$. All of these
algorithms compute the roots accurately enough to separate them.

\paragraph{Truncation.}
To make the algorithm efficient when $b<n$, we truncate $f$ to its top
$\min(n,b)+1$ coefficients (Step~1). By Lemma~\ref{lem:leading-part},
every integer root $r$ of $f$ with $|r|\geq2$ lies within distance $2$ of
a root of the truncated polynomial $g$, so it suffices to approximate the
roots of $g$. This costs $\softO(d(b+d))$ bit operations, instead of
$\softO(n(n+b))$ for $f$ itself. Truncation has been used
before~\cite{CKS99,Mor22}: Cucker, Koiran, and Smale split off the
low-order terms of a sparse polynomial across a large gap in the
exponents to find its integer roots, and Moroz approximates a polynomial
on small disks by truncated expansions to speed up root isolation and
multipoint evaluation.

The rest of the paper is organized as follows.
Section~\ref{sec:estimates} gives the bounds on polynomial roots needed for the algorithm. Section~\ref{sec:algorithm} fills in the
details of each step and proves part~(b) of
Theorem~\ref{thm:integer-roots}, and Appendix~\ref{sec:elementary}
proves part~(a).

\section{Bounds on roots of polynomials}
\label{sec:estimates}

We define the relevant polynomial notation and derive bounds on the roots of truncated and perturbed polynomials. For a polynomial $p(x)=\sum_jp_jx^j\in\CC[x]$, we use the standard
coefficient norms $\|p\|_1=\sum_j|p_j|$,
$\|p\|_2=\bigl(\sum_j|p_j|^2\bigr)^{1/2}$, and
$\|p\|_\infty=\max_j|p_j|$. To bound the sum of the bit lengths of the
roots of a polynomial, we recall its Mahler measure. For a nonzero
polynomial $p(x)=c\prod_{i=1}^d(x-\beta_i)\in\CC[x]$, the \emph{Mahler
measure} of $p$ is defined as
\[
 M(p)=|c|\prod_{i=1}^d\max(1,|\beta_i|).
\]
We use the following inequality of Landau~\cite{Landau1905},
which bounds the Mahler measure by the coefficient $\ell_2$ norm.

\begin{lemma}\label{lem:landau}
For every nonzero polynomial $p\in\CC[x]$, we have $M(p)\leq\|p\|_2$.
\end{lemma}

\begin{proof}
Write $p(x)=\sum_{j=0}^{m} p_j x^j$.
By Jensen's formula and the concavity of the logarithm,
\[
\log M(p)
=
\frac{1}{2\pi}\int_0^{2\pi}\log |p(e^{\iu\theta})|\,d\theta
\leq
\frac12\log\left(
\frac{1}{2\pi}\int_0^{2\pi}|p(e^{\iu\theta})|^2\,d\theta
\right).
\]
By Parseval's identity,
\[
\frac{1}{2\pi}\int_0^{2\pi}|p(e^{\iu\theta})|^2\,d\theta
=
\sum_{j=0}^{m}|p_j|^2
=
\|p\|_2^2.
\]
Here the first equality follows from the orthogonality of
the functions $e^{\iu j\theta}$ on $[0,2\pi]$.
Thus $\log M(p)\leq \log\|p\|_2$, proving the claim.
\end{proof}

Next we show that any integer root of $f$ outside $\{-1,0,1\}$
must be close to a complex root of the truncated polynomial $g$.

\begin{lemma}\label{lem:leading-part}
Let $d=\min(n,b)$ and $g(x)=a_nx^d+a_{n-1}x^{d-1}+\cdots+a_{n-d}$. If
$d<n$, then for every integer root $r$ of $f$ with $|r|\geq2$, there is
a complex root $\alpha$ of $g$ such that $|r-\alpha|<2$.
\end{lemma}

\begin{proof}
Since $d<n$, we have $d=b$. Write $f(x)=x^{n-d}g(x)+h(x)$ with
$h(x)=\sum_{k=0}^{n-d-1}a_kx^k$, and let $r$ be an integer root of $f$
with $|r|\geq2$. Because $f(r)=0$ and $|a_k|<2^b$, we have
$g(r)=-h(r)/r^{n-d}$, which gives
\[
 |g(r)|=\frac{|h(r)|}{|r|^{n-d}}
 <\frac{2^b}{|r|^{n-d}}\sum_{k=0}^{n-d-1}|r|^k
 <\frac{2^b}{|r|-1}\leq2^b.
\]
Now let $\alpha_1,\ldots,\alpha_d$ be the complex roots of $g$, counted
with multiplicity. If $|r-\alpha_j|\geq2$ held for all $j$, then since
$a_n$ is a nonzero integer we would have
\[
 |g(r)|=|a_n|\prod_{j=1}^d|r-\alpha_j|\geq2^d=2^b,
\]
contradicting $|g(r)|<2^b$.
\end{proof}

The assumption $|r|\geq2$ in Lemma~\ref{lem:leading-part} cannot be
dropped. For example, if $f(x)=x^5+3x^4+3x^3-3x^2-3x-1$ and $b=2$, the
root $r=1$ has distance $\sqrt7>2$ from both roots of the truncated
polynomial $g(x)=x^2+3x+3$. To take care of small roots, we simply put
$0$ and $\pm1$ into $\mathcal C$ by hand.

Next, we consider a polynomial $\widehat g\in\CC[x]$ whose coefficients
are close to those of $g$ in the $\ell_1$ norm. The following lemma
shows that every root of $g$ lies near some root of $\widehat g$, and at
the same time bounds the sum of the logarithmic sizes of the roots of
$\widehat g$.

\begin{lemma}\label{lem:perturbation}
Let $g\in\ZZ[x]$ be a polynomial of degree $d\geq1$ with
$\|g\|_\infty<2^b$. Let $\widehat g\in\CC[x]$ be a polynomial of degree
$d$ with roots $\beta_1,\ldots,\beta_d$, counted with multiplicity. If
\[
 \delta:=\|\widehat g-g\|_1<\frac1{8\cdot9^d(2^b\sqrt{d+1}+1)},
\]
then every root of $g$ lies within distance $1/8$ of at least one
$\beta_i$, and the roots $\beta_i$ satisfy
\[
 \sum_{i=1}^d\log_2\max(1,|\beta_i|)
 \leq b+\tfrac12\log_2(d+1)+2.
\]
\end{lemma}

\begin{proof}
Let $c$ be the leading coefficient of $\widehat g$. Since the leading
coefficient of $g$ is a nonzero integer and $\delta<1/2$, we have
$|c|>1/2$. By Lemma~\ref{lem:landau},
\begin{equation}\label{eq:mahler-bound}
 |c|\prod_{i=1}^d\max(1,|\beta_i|)=M(\widehat g)
 \leq\|\widehat g\|_2\leq\|g\|_2+\delta
 \leq2^b\sqrt{d+1}+\delta.
\end{equation}
The second claim follows right away from~\eqref{eq:mahler-bound}: the
sum on the left-hand side is $\log_2(M(\widehat g)/|c|)$, and since
$|c|>1/2$ and $\delta<1$, this is bounded above by
\[
 \log_2\bigl(2(2^b\sqrt{d+1}+1)\bigr)\leq b+\tfrac12\log_2(d+1)+2.
\]

For the first claim, let $\alpha$ be a root of $g$ and put
$A=\max(1,|\alpha|)$. Suppose for contradiction that
$|\alpha-\beta_i|\geq1/8$ for all $i=1,\ldots,d$. Then for each $i$,
\[
 A\leq\max(1,|\beta_i|)\,(1+|\alpha-\beta_i|)
 \leq9\,|\alpha-\beta_i|\max(1,|\beta_i|).
\]
Multiplying these inequalities over all $i$ and
using~\eqref{eq:mahler-bound}, we find
\[
 |\widehat g(\alpha)|
 =|c|\prod_{i=1}^d|\alpha-\beta_i|
 \geq\frac{|c|A^d}{9^d\prod_{i=1}^d\max(1,|\beta_i|)}
 =\frac{|c|^2A^d}{9^dM(\widehat g)}
 \geq\frac{A^d}{4\cdot9^d(2^b\sqrt{d+1}+\delta)}.
\]
On the other hand, since $g(\alpha)=0$, we have
$|\widehat g(\alpha)|=|(\widehat g-g)(\alpha)|\leq\delta A^d$. Comparing
the upper and lower bounds on $|\widehat g(\alpha)|$ gives
$\delta(2^b\sqrt{d+1}+\delta)\geq1/(4\cdot9^d)$, which contradicts our
hypothesis on $\delta$.
\end{proof}

\section{The algorithm}
\label{sec:algorithm}
\subsection{Algorithmic tools}\label{sec:known}

Our algorithms rely on two kinds of tools. The first are standard
algorithms for exact polynomial arithmetic. The second is Pan's numerical
algorithm for approximate factorization, which we need in a specific
form and therefore state separately as Theorem~\ref{thm:pan}.

\paragraph{Exact polynomial arithmetic.}
We use multipoint evaluation, evaluation at a single integer, and Taylor
shifts. We list them below, with precise references to the textbook of
von zur Gathen and Gerhard~\cite{vzGG13} where applicable.

\begin{enumerate}
    \item\label{comp:multipoint_evaluation} A polynomial $f$ of degree $<n$ over a ring can be
    evaluated at $m$ points $a_1,\ldots,a_m$ in $\softO(n+m)$ ring operations \cite[Corollary~10.8]{vzGG13}.

    \item\label{comp:substitute_in_poly} Let $f \in \mathbb{Z}[x]$ have degree $n$ with $b$-bit
    coefficients, and let $x \in \mathbb{Z}$ have $b$ bits. Then $f(x)$ can be computed in
    $\softO(nb)$ bit operations: writing $f = f_0 + x^{\lceil n/2 \rceil} f_1$
    and recursing yields $T(n) = 2T(n/2) + \softO(nb) = \softO(nb)$.

    \item\label{comp:taylor_shift} Let $R$ be a commutative ring, $f \in R[x]$ of degree $\le n$,
    and $a,b \in R$. Then $f(ax+b)$ can be computed in $\softO(n)$ ring operations: the polynomial
    $g(x)=f(x+b)$ can be computed in $\softO(n)$ ring operations \cite[Corollary~9.16]{vzGG13},
    and $f(ax+b)=g(ax)$ is obtained from $g$ by multiplying its $k$th coefficient by $a^k$.
\end{enumerate}

As a corollary of (1) and (3), we can evaluate and shift polynomials
over $\ZZ$ or $\ZZ[\iu]$ exactly, by working in $\ZZ/2^L\ZZ$ or
$(\ZZ/2^L\ZZ)[\iu]$ for a suitable $L$. A ring operation there takes
$\softO(L)$ bit operations, and an element of $\ZZ$ or $\ZZ[\iu]$ whose real and imaginary parts have absolute value less than
$2^{L-1}$ is determined by its residue modulo $2^L$: it is obtained by
taking the residues of smallest absolute value. Hence, if the real and
imaginary parts of all outputs have absolute value less than
$2^{L-1}$, then (1) evaluates a polynomial of degree $<n$ over
$\ZZ[\iu]$ at $m$ points in $\softO((n+m)L)$ bit operations, and (3)
computes $f(x+u)$ for $f\in\ZZ[\iu][x]$ of degree $\le n$ and
$u\in\ZZ[\iu]$ in
$\softO(nL)$ bit operations.

Throughout, we refer to the items of the list above as \cref{sec:known}(1)--(3).

\paragraph{Approximate factorization.}
We call a complex number of the form $(u+\iu v)/2^k$ with $u,v\in\ZZ$
and $k\geq0$ a \emph{dyadic Gaussian rational}. The \emph{bit length} of
an integer $u$ is $\lceil\log_2(|u|+1)\rceil$, and the \emph{bit size}
of $(u+\iu v)/2^k$ is the sum of the bit lengths of $u$, $v$, and $k$.
More generally, the bit size of $(u+\iu v)/w$ with $u,v,w\in\ZZ$ and
$w\geq1$ is the sum of the bit lengths of $u$, $v$, and $w$.

To approximate the roots of $g$, we rely on Pan's approximate
factorization algorithm~\cite[Theorem~2.1.1]{Pan02}. We use the version
stated by Mehlhorn, Sagraloff, and Wang~\cite[Theorem~2]{MSW15}, which
gives explicit bit complexity bounds. Specialized to exact dyadic inputs,
their result reads as follows.

\begin{theorem}[{\cite[Theorem~2]{MSW15}}]\label{thm:pan}
Let $p$ be a polynomial of degree $d\geq1$ whose roots all lie in the
closed unit disk, and whose coefficients are dyadic Gaussian rationals
with a common denominator $2^s$, where $s=O(P)$, and numerators of
$O(P)$ bits, for an integer $P\geq d\log_2(d+1)$. Then one can
deterministically compute, in $\softO(dP)$ bit operations, dyadic
Gaussian rationals $w_1,\ldots,w_d$ with a common denominator $2^B$,
where $B=O(P)$, such that
\[
 \Bigl\|p-p_d\prod_{j=1}^d(x-w_j)\Bigr\|_1\leq2^{-P}\|p\|_1.
\]
\end{theorem}

\begin{remark}\label{rem:pan}
Theorem~\ref{thm:pan} is only used in the proof of
Theorem~\ref{thm:integer-roots}(b). While preparing this paper, we were
unable to verify parts of the running time analysis
in~\cite{Pan02}. In Stage~1 of~\cite[Algorithm~2.6.1]{Pan02}, the
algorithm calls itself on a higher-order derivative of its input with
the precision parameter $B$ replaced by $2B\log s$, and we could not
verify that the resulting growth of the working precision over the
recursion, together with the cost of the $O(\log n)$ passes through
Stage~3, stays within the bound of~\cite[Theorem~2.1.1]{Pan02}.
Theorem~\ref{thm:integer-roots}(a) does not
depend on any of this: its proof in Appendix~\ref{sec:elementary} uses
only Lemmas~\ref{lem:landau}, \ref{lem:leading-part},
and~\ref{lem:multipoint} from the main text.
\end{remark}

\subsection{Proof of Theorem~\ref{thm:integer-roots}(b)}

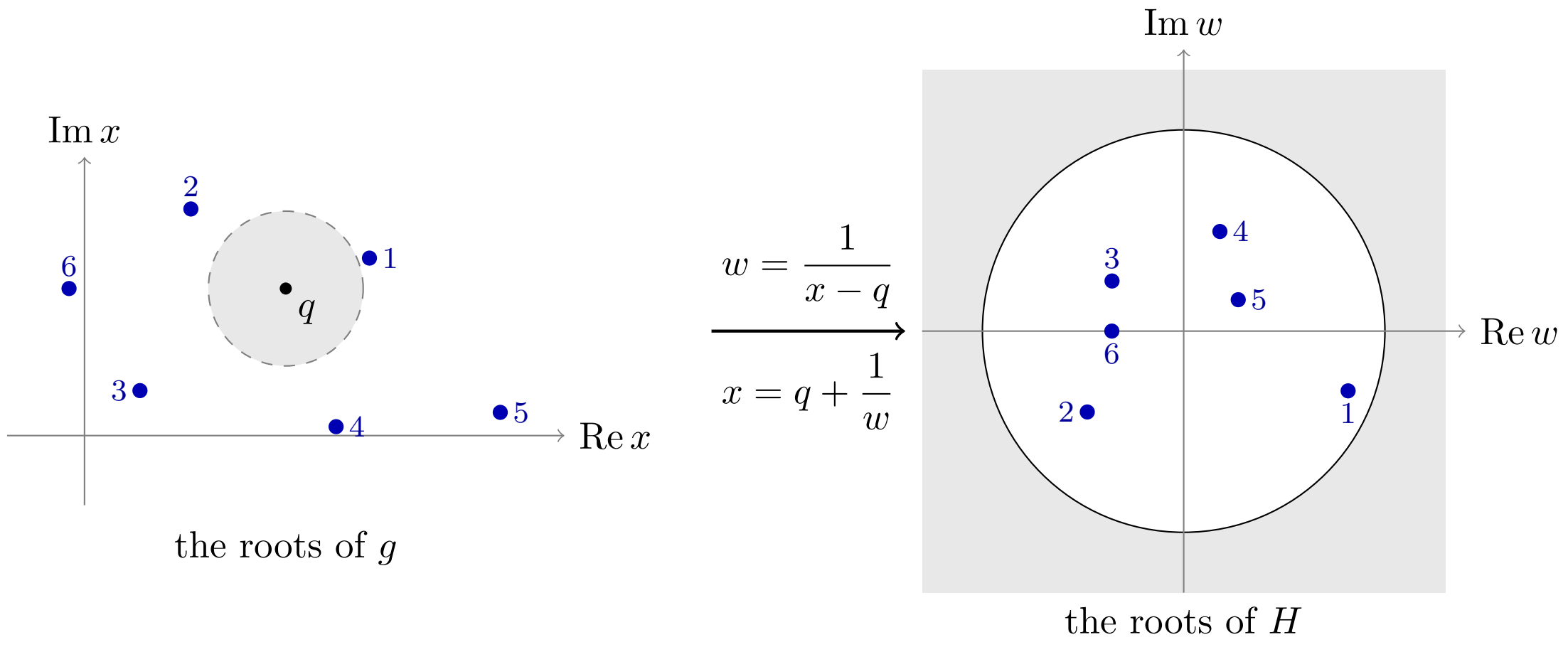
\begin{figure}[ht]
\centering
\begin{tikzpicture}[x=0.72cm,y=0.72cm,font=\small,
  rt/.style={circle,fill=blue!70!black,inner sep=0pt,minimum size=4pt},
  lb/.style={font=\scriptsize,text=blue!60!black,inner sep=3.2pt}]
\begin{scope}
\fill[gray!18] (2.600,1.900) circle (1);
\draw[dashed,gray] (2.600,1.900) circle (1);
\draw[->,gray] (-1.00,0)--(6.20,0) node[right,black]{$\operatorname{Re}x$};
\draw[->,gray] (0,-0.90)--(0,3.60) node[above,black]{$\operatorname{Im}x$};
\fill (2.600,1.900) circle (1.6pt) node[below right=-1pt]{$q$};
\node[rt] at (3.681,2.293) {}; \node[lb,right] at (3.681,2.293) {1};
\node[rt] at (1.374,2.928) {}; \node[lb,above] at (1.374,2.928) {2};
\node[rt] at (0.716,0.581) {}; \node[lb,left] at (0.716,0.581) {3};
\node[rt] at (3.250,0.115) {}; \node[lb,right] at (3.250,0.115) {4};
\node[rt] at (5.371,0.300) {}; \node[lb,right] at (5.371,0.300) {5};
\node[rt] at (-0.200,1.900) {}; \node[lb,above] at (-0.200,1.900) {6};
\node[anchor=north] at (2.60,-1.05) {the roots of $g$};
\end{scope}
\begin{scope}[shift={(14.200,1.350)},x=1.8719999999999999cm,y=1.8719999999999999cm]
\fill[gray!18] (-1.3,-1.3) rectangle (1.3,1.3);
\fill[white] (0,0) circle (1); \draw (0,0) circle (1);
\draw[->,gray] (-1.3,0)--(1.4,0) node[right,black]{$\operatorname{Re}w$};
\draw[->,gray] (0,-1.3)--(0,1.4) node[above,black]{$\operatorname{Im}w$};
\node[rt] at (0.817,-0.297) {}; \node[lb,below] at (0.817,-0.297) {1};
\node[rt] at (-0.479,-0.402) {}; \node[lb,left] at (-0.479,-0.402) {2};
\node[rt] at (-0.356,0.249) {}; \node[lb,above] at (-0.356,0.249) {3};
\node[rt] at (0.180,0.495) {}; \node[lb,right] at (0.180,0.495) {4};
\node[rt] at (0.271,0.156) {}; \node[lb,right] at (0.271,0.156) {5};
\node[rt] at (-0.357,-0.000) {}; \node[lb,below] at (-0.357,-0.000) {6};
\node[anchor=north] at (0,-1.3) {the roots of $H$};
\end{scope}
\draw[->,thick] (8.10,1.35) -- node[above=2pt]{$w=\dfrac{1}{x-q}$} node[below=2pt]{$x=q+\dfrac1w$} (10.60,1.35);
\end{tikzpicture}
\caption{The substitution $x=q+1/w$. All roots of $g$ lie in the bright
area $|x-q|\geq1$, which is mapped onto the bright disk $|w|\leq1$.
Matching numbers mark corresponding roots.}
\label{fig:inversion}
\end{figure}

Before we can apply Theorem~\ref{thm:pan}, we need to map the roots of
$g$ into the unit disk. For this, it suffices to find a Gaussian integer
$q$ whose distance to every root of $g$ is greater than $1$, so that the
polynomial $w^dg(q+1/w)$ has all its roots inside the unit disk
(see \cref{fig:inversion}). The next lemma shows how to find such a $q$.

\begin{lemma}\label{lem:shift}
Let $g(x)=\sum_{m=0}^dg_mx^m\in\ZZ[x]$ be a polynomial of degree
$d\geq1$ with $|g_m|<2^b$. In $\softO(d(b+d))$ bit operations, we can
deterministically find a Gaussian integer $q$ with $|q|\leq10d^2+1$
whose distance to every complex root of $g$ is strictly greater than
$1$, and compute the coefficients of the shifted polynomial $g(q+x)$,
each of which has bit length $\softO(b+d)$.
\end{lemma}

\begin{proof}
Let $R=10d^2$, and let $\zeta_k=e^{2\pi\iu k/(d+1)}$ for
$0\leq k\leq d$. The computation goes as follows.
\begin{enumerate}
 \item Compute each $R\zeta_k$ to within $1/8$, and round its real and
 imaginary parts to the nearest integers. This gives Gaussian integers
 $q_0,\ldots,q_d$ with $|q_k-R\zeta_k|\leq1$.
 \item Compute $g(q_0),\ldots,g(q_d)$, and let $q$ be a point $q_k$ for
 which $|g(q_k)|$ is largest.
 \item Compute the coefficients of the shifted polynomial $g(q+x)$.
\end{enumerate}

First, we prove that $g$ has no complex root $\alpha$ with
$|\alpha-q|\leq1$. Define
\[
 A_R=\Bigl(\sum_{m=0}^d|g_m|^2R^{2m}\Bigr)^{1/2}.
\]
Because $\deg g<d+1$, the average of $|g(R\zeta_k)|^2$ over
$0\leq k\leq d$ is $A_R^2$. Indeed, expanding the square and using
$\overline{\zeta_k}=\zeta_k^{-1}$, we get
\[
 \frac1{d+1}\sum_{k=0}^d|g(R\zeta_k)|^2
 =\sum_{m,l=0}^dg_m\overline{g_l}\,R^{m+l}\cdot\frac1{d+1}\sum_{k=0}^d\zeta_k^{m-l}
 =\sum_{m=0}^d|g_m|^2R^{2m}=A_R^2,
\]
since for $0\leq m,l\leq d$ the inner average is a geometric series with
ratio $e^{2\pi\iu(m-l)/(d+1)}$, which equals $1$ if $m=l$ and $0$
otherwise. Hence there is at least one index $k_0$ for which
$|g(R\zeta_{k_0})|\geq A_R$.

Next, we bound $|g'(z)|$ from above in the disk $|z|\leq R+2$. For such
$z$,
\begin{align*}
 |g'(z)|
 &=\Bigl|\sum_{m=1}^dm\,g_mz^{m-1}\Bigr|
 \leq d\sum_{m=1}^d|g_m|(R+2)^{m-1}
 \leq\frac dR\Bigl(1+\frac2R\Bigr)^d\sum_{m=1}^d|g_m|R^m\\
 &\leq\frac{d^{3/2}}{R}\Bigl(1+\frac2R\Bigr)^{R/2}A_R
 \leq\frac{e\,d^{3/2}}{R}A_R
 \leq\frac{A_R}{3}.
\end{align*}
Here the third inequality uses $d\leq R/2$ and the Cauchy--Schwarz
inequality $\sum_{m=1}^d|g_m|R^m\leq\sqrt d\,A_R$, and the last one uses
$R=10d^2$. Hence, for any $z,w$ with $|z|,|w|\leq R+2$, we have
\[
 |g(z)-g(w)|=\Bigl|\int_w^zg'(t)\,dt\Bigr|\leq\frac{A_R}{3}|z-w|.
\]
We apply this twice. First, $|R\zeta_{k_0}|=R$ and $|q_{k_0}|\leq R+1$,
so $|g(q_{k_0})|\geq|g(R\zeta_{k_0})|-A_R/3\geq2A_R/3$, and by the
choice of $q$ also $|g(q)|\geq2A_R/3$. Second, if $g$ had a root $\alpha$
with $|\alpha-q|\leq1$, then $|\alpha|\leq R+2$ and
$|g(q)|=|g(q)-g(\alpha)|\leq A_R/3$, a contradiction.

It remains to bound the running time. Step~1 needs $\cos$, $\sin$ and
$\pi$ to $O(\log d)$ bits of precision, which takes
$\operatorname{poly}(\log d)$ bit operations per
point~\cite[Sections~4.9.1 and~4.12]{brent2010modern}, so $\softO(d)$
in total. Let $K$ be the smallest integer with
$2^{K-1}>2^b(d+1)(R+2)^d$, so that $K=\softO(b+d)$. Since
$|q_k|\leq R+1$, we have
$|g(q_k)|\leq\sum_{m=0}^d|g_m||q_k|^m<2^b(d+1)(R+1)^d<2^{K-1}$. Hence \tool{comp:multipoint_evaluation}, applied over
$(\ZZ/2^K\ZZ)[\iu]$, carries out Step~2 in $\softO(dK)$ bit operations.
For Step~3 we use \tool{comp:taylor_shift}, again over
$(\ZZ/2^K\ZZ)[\iu]$. The coefficient of $x^j$ in $g(q+x)$ is
$\sum_{m\geq j}\binom mjg_mq^{m-j}$, whose absolute value is at most
$2^b\sum_{m=0}^d(1+|q|)^m\leq2^b(d+1)(R+2)^d<2^{K-1}$, so Step~3 also
takes $\softO(dK)=\softO(d(b+d))$ bit operations, and each coefficient
of $g(q+x)$ has bit length at most $K$.
\end{proof}

Putting Lemma~\ref{lem:shift}, Theorem~\ref{thm:pan}, and
Lemma~\ref{lem:perturbation} together gives the root approximation
procedure needed for Step~2.

\begin{proposition}\label{prop:approximate-roots}
Let $g\in\ZZ[x]$ be a polynomial of degree $d\geq1$ with coefficients of
absolute value less than $2^b$. In $\softO(d(b+d))$ bit operations, one
can deterministically compute complex numbers $\beta_1,\ldots,\beta_d$
with rational real and imaginary parts, of total bit size
$\softO(d(b+d))$, such that
every root of $g$ has distance less than $1/8$ to at least one
$\beta_i$, and
\[
 \sum_{i=1}^d\log_2\max(1,|\beta_i|)\leq b+\tfrac12\log_2(d+1)+2.
\]
\end{proposition}

\begin{proof}
Let $R=10d^2$, and let $q$ and $g(q+x)$ be as computed by
Lemma~\ref{lem:shift}, so that $|q|\leq R+1$. The computation goes as
follows.
\begin{enumerate}
 \item Form the polynomial $H(w)=w^dg(q+1/w)$. If
 $g(q+x)=\sum_{j=0}^dc_jx^j$, then $H(w)=\sum_{j=0}^dc_jw^{d-j}$, so the
 coefficients of $H$ are those of $g(q+x)$ in reverse order.
 \item Apply Theorem~\ref{thm:pan} to $H$ with precision
 \[
  P=2b+8d+16+\bigl\lceil2d\log_2(R+2)+\log_2(d+1)\bigr\rceil,
 \]
 and obtain $w_1,\ldots,w_d$.
 \item Output $\beta_i=q+1/w_i$ for $1\leq i\leq d$.
\end{enumerate}

First, $H$ satisfies the hypotheses of Theorem~\ref{thm:pan}. Since $g$
has no root within distance $1$ of $q$, the polynomial $H$ has degree
$d$, and its roots $1/(\alpha-q)$, where $\alpha$ runs over the roots
of $g$, lie in the open unit disk. Its coefficients are Gaussian
integers, and since $\|(x+q)^m\|_1=(1+|q|)^m$, we have
\[
 \|H\|_1=\|g(q+x)\|_1
 \leq2^b\sum_{m=0}^d(1+|q|)^m
 \leq(d+1)\,2^b(R+2)^d
 \leq2^P,
\]
so they have $O(P)$ bits. Finally, $P\geq2d\log_2(R+2)\geq d\log_2(d+1)$.

To go back from $w$ to $x$, we use the map $E\mapsto E^*$, where
$E^*(x)=(x-q)^dE(1/(x-q))$ for a polynomial $E$ of degree at most $d$.
This map is linear and sends $H$ to $g$. Moreover, it sends
$E(w)=\sum_je_jw^j$ to $\sum_je_j(x-q)^{d-j}$, and
$\|(x-q)^k\|_1=(1+|q|)^k$, so
$\|E^*\|_1\leq(1+|q|)^d\|E\|_1\leq(R+2)^d\|E\|_1$.

Next, we prove that the $\beta_i$ have the required properties. Let
$H_d=g(q)\neq0$ be the leading coefficient of $H$, let
$\widehat H(w)=H_d\prod_{i=1}^d(w-w_i)$, so that
$\|\widehat H-H\|_1\leq2^{-P}\|H\|_1$ by Theorem~\ref{thm:pan}, and let
\[
 \widehat g(x)=\widehat H^*(x)=H_d\prod_{i=1}^d\bigl(1-w_i(x-q)\bigr).
\]
By Lemma~\ref{lem:perturbation}, applied to $g$ and $\widehat g$, it
suffices to show that
$\|\widehat g-g\|_1<1/\bigl(8\cdot9^d(2^b\sqrt{d+1}+1)\bigr)$. Indeed,
\begin{align*}
 \|\widehat g-g\|_1=\|(\widehat H-H)^*\|_1
 &\leq(R+2)^d\|\widehat H-H\|_1
 \leq2^{-P}(R+2)^d\|H\|_1\\
 &\leq2^{-P}(d+1)\,2^b(R+2)^{2d}
 \leq2^{-b-8d-16}\\
 &<\frac1{2^3\cdot2^{4d}\cdot2^{b+d+1}}
 \leq\frac1{8\cdot9^d(2^b\sqrt{d+1}+1)}.
\end{align*}
In particular, the coefficients of $x^d$ in $\widehat g$ and $g$ differ
by less than $1/2$. The first is $H_d\prod_{i=1}^d(-w_i)$ and the second
is a nonzero integer, so every $w_i$ is nonzero, and $\widehat g$ has
degree $d$ and roots $\beta_i=q+1/w_i$. This proves that every root of $g$ has distance less than $1/8$ to some
$\beta_i$, and the bound on $\sum_i\log_2\max(1,|\beta_i|)$.

Finally, computing $q$ and $g(q+x)$ takes $\softO(d(b+d))$ bit
operations by Lemma~\ref{lem:shift}, Step~1 only reverses the order of
the coefficients, and Step~2 takes $\softO(dP)=\softO(d(b+d))$ bit
operations by Theorem~\ref{thm:pan}. The $w_i$ are written down within
this time, so they have total bit size $\softO(d(b+d))$, and computing
the $\beta_i$ exactly in Step~3 takes $\softO(d(b+d))$ bit operations
and gives numbers of total bit size $\softO(d(b+d))$.
\end{proof}

We now turn to the exact multipoint evaluation needed in Step~4. The key
feature of the lemma below is that the running time depends on the sum
of the bit lengths of the evaluation points rather than on the maximum
bit length. Combined with the bound $S=O(b)$, this avoids the
verification bottleneck discussed in~\cite[Appendix~A]{HH22}.

\begin{lemma}\label{lem:multipoint}
Let $f\in\ZZ[x]$ have degree $n$ and coefficients of absolute value less
than $2^b$, let $c_1,\ldots,c_m\in\ZZ$ with $m=O(n)$, and let
$S=\sum_{i=1}^m\lceil\log_2(2+|c_i|)\rceil$. Then the values
$f(c_1),\ldots,f(c_m)$ can be computed deterministically in
$\softO(nb+nS)$ bit operations.
\end{lemma}

\begin{proof}
By splitting the points into $O(1)$ groups of at most $n$ points, we
may assume $m\leq n$. We sort the points into buckets according to
$\lceil\log_2(2+|c_i|)\rceil$: for $\tau=1,2,4,\ldots$ up to the first power of two that is
at least $S$, the bucket of $\tau$ contains the points $c_i$ with
$\tau/2<\lceil\log_2(2+|c_i|)\rceil\leq\tau$. Every point lies in
exactly one of these $O(\log S)$ buckets. If the bucket of $\tau$
contains $t$ points, then each of them satisfies $|c_i|<2^\tau$, and
$t\tau/2<S$. Hence it is enough to show how to evaluate $f$ at $t\leq n$
points satisfying $|c_i|<2^\tau$ in $\softO(nb+nt\tau)$ bit operations:
summing over the buckets then gives $\softO(nb+nS)$, where the $\softO$
notation absorbs the $O(\log S)$ factor coming from the number of
buckets.

Break $f$ into blocks of $t$ coefficients,
\[
 f(x)=\sum_{j=0}^{\lfloor n/t\rfloor}x^{tj}F_j(x),
\]
where $\deg F_j<t$. Since $t\leq n$, there are at most $2n/t$ blocks.
Each $F_j$ has $t$ coefficients of absolute value less than $2^b$, so
$|F_j(c_i)|<t\,2^b2^{t\tau}\leq2^{L-1}$, where
$L=b+t\tau+\lceil\log_2(t+1)\rceil+3$. Hence, for each $j$,
\tool{comp:multipoint_evaluation}, applied over $\ZZ/2^L\ZZ$, evaluates
$F_j$ at all $t$ points in $\softO(tL)$ bit operations, which is
$\softO(nL)=\softO(n(b+t\tau))$ over all blocks.

To put the block values together at a point $c_i$, we compute $c_i^t$ by
repeated squaring in $\softO(t\tau)$ bit operations. Then
$f(c_i)=\sum_jF_j(c_i)\,(c_i^t)^j$ is the value at $c_i^t$, which has at
most $t\tau\leq L$ bits, of a polynomial with at most $2n/t$
coefficients of at most $L$ bits each. By
\tool{comp:substitute_in_poly}, this takes
$\softO(nL/t)=\softO(nb/t+n\tau)$ bit operations for a single point,
and therefore $\softO(nb+nt\tau)$ bit operations for all $t$ points.
\end{proof}

\begin{proof}[Proof of Theorem~\ref{thm:integer-roots}(b)]
We follow Steps~1--4 of the algorithm. Recall that $d=\min(n,b)$ and
$g(x)=a_nx^d+a_{n-1}x^{d-1}+\cdots+a_{n-d}$ (Step~1). Since $g$ has
integer coefficients of absolute value less than $2^b$, we can apply
Proposition~\ref{prop:approximate-roots} to $g$ and obtain
$\beta_1,\ldots,\beta_d$ (Step~2). For each $i$, let $z_i$ be an integer
nearest to $\operatorname{Re}\beta_i$, and form the candidate set
\[
 \mathcal C=\{0,\pm1\}\cup\{z_i+j : 1\leq i\leq d,\ |j|\leq3\}
\]
(Step~3). Finally, we compute $f(c)$ for all $c\in\mathcal C$ using
Lemma~\ref{lem:multipoint}, and output those $c$ with $f(c)=0$
(Step~4).

To see that the algorithm is correct, note first that Step~4 checks
$f(c)=0$ exactly, so only true roots are output. Conversely, let
$r\in\ZZ$ be an integer root of $f$. If $|r|\leq1$, then $r\in\mathcal
C$ by definition. If $|r|\geq2$, then $g$ has a complex root $\alpha$
with $|r-\alpha|<2$ (by Lemma~\ref{lem:leading-part} when $d<n$, and
trivially when $d=n$ since $g=f$). By
Proposition~\ref{prop:approximate-roots}, there is some $\beta_i$ with
$|\alpha-\beta_i|<1/8$. Since $r$ is real and $z_i$ is an integer
closest to $\operatorname{Re}\beta_i$, we get
\[
 |r-z_i|\leq|r-\operatorname{Re}\beta_i|+\tfrac12
 \leq|r-\beta_i|+\tfrac12
 \leq|r-\alpha|+|\alpha-\beta_i|+\tfrac12
 <2+\tfrac18+\tfrac12<3.
\]
Since $r-z_i$ is an integer, it follows that $r\in\mathcal C$, so $r$ is
output in Step~4.

Next, let us bound the total bit length $S$ of the candidate set
$\mathcal C$. We have $|\mathcal C|\leq7d+3$. Every candidate $c=z_i+j$
with $|j|\leq3$ satisfies $|c|\leq|\beta_i|+7/2$, which implies
$\lceil\log_2(2+|c|)\rceil\leq4+\log_2\max(1,|\beta_i|)$. Using
Proposition~\ref{prop:approximate-roots} and the fact that $d\leq b$, we
obtain
\[
 S\leq6+7\sum_{i=1}^d\bigl(4+\log_2\max(1,|\beta_i|)\bigr)
 \leq7b+28d+\tfrac72\log_2(d+1)+20=O(b).
\]

Finally, we add up the running times of the steps. Steps~1 and~3 only
write down $g$ and round the $\beta_i$, which takes $\softO(nb)$ bit
operations. By Proposition~\ref{prop:approximate-roots}, Step~2 runs in
$\softO(d(b+d))\subseteq\softO(nb)$ bit operations. Since
$|\mathcal C|\leq7d+3=O(n)$ and $S=O(b)$, Lemma~\ref{lem:multipoint}
carries out Step~4 in $\softO(nb+nS)=\softO(nb)$ bit operations.
\end{proof}

\appendix
\section{An elementary candidate search}
\label{sec:elementary}
In this appendix we prove Theorem~\ref{thm:integer-roots}(a). We keep
the notation of the paper: $d=\min(n,b)$, and
$g(x)=\sum_{j=0}^dg_jx^j=a_nx^d+\cdots+a_{n-d}$, so that $|g_j|<2^b$ and
$|g_d|\geq1$. Instead of approximating the complex roots of $g$, we
search for candidate integers directly, by bisection on the real line.
To discard an interval, we use the exclusion test of Becker, Sagraloff,
Sharma, and Yap, which combines the case $k=0$ of Pellet's theorem with
Graeffe's root squaring~\cite[Section~3]{BSSY18}. We include the short
proofs for completeness. From the rest of the paper, this appendix uses
only Lemmas~\ref{lem:landau}, \ref{lem:leading-part},
and~\ref{lem:multipoint}, and the Taylor shift
(\tool{comp:taylor_shift}). In particular, it does not use
Theorem~\ref{thm:pan}.

\begin{lemma}\label{lem:el-test}
Let $q=\sum_{k=0}^dq_kx^k\in\CC[x]$ have degree $d\geq1$.
\begin{enumerate}[\rm(a)]
 \item If $|q_0|>\sum_{k=1}^d|q_k|$, then $q$ has no root in the disk
 $|x|\leq1$.
 \item Otherwise, $q$ has a root of absolute value at most $1.5d$.
\end{enumerate}
\end{lemma}

\begin{proof}
(a) For $|x|\leq1$ we have
$|q(x)|\geq|q_0|-\sum_{k\geq1}|q_k||x|^k>0$.

(b) Suppose that all roots $\gamma_1,\ldots,\gamma_d$ of $q$ satisfy
$|\gamma_i|>1.5d$. Then $q_0\neq0$ and
$q(x)=q_0\prod_{i=1}^d(1-x/\gamma_i)$. Up to sign, the coefficient of
$x^k$ in the product is the sum of the $\binom dk$ products of $k$ of the
numbers $1/\gamma_i$, so $|q_k|\leq|q_0|\binom dk(1.5d)^{-k}$. Hence
\[
 \sum_{k=1}^d|q_k|\leq|q_0|\Bigl(\bigl(1+\tfrac1{1.5d}\bigr)^d-1\Bigr)
 \leq|q_0|\bigl(e^{2/3}-1\bigr)<|q_0|,
\]
and we are in case (a).
\end{proof}

To shrink the factor $1.5d$ in Lemma~\ref{lem:el-test}, we square the
roots a few times. If $p(x)=c\prod_{i=1}^d(x-\delta_i)$, then
$p(x)p(-x)=(-1)^dc^2\prod_{i=1}^d(x^2-\delta_i^2)$, so there is a
polynomial $Gp$ with $Gp(x^2)=(-1)^dp(x)p(-x)$, namely
$Gp(y)=c^2\prod_{i=1}^d(y-\delta_i^2)$. Its roots are the squares of the
roots of $p$.

Let $N=\lceil\log_2\log_2(2d)\rceil$, so that $2^N\geq\log_2(2d)$. For
integers $m$ and $r\geq1$, the test $\mathrm{Test}(m,r)$ computes
$p(x)=g(m+rx)$ and $q=G^Np$, and answers \emph{empty} if
$|q_0|>\sum_{k\geq1}|q_k|$.

\begin{lemma}\label{lem:el-soft}
If $\mathrm{Test}(m,r)$ answers empty, then every root $\alpha$ of $g$
satisfies $|\alpha-m|>r$. Otherwise, some root $\alpha$ of $g$ satisfies
$|\alpha-m|\leq2r$. If $|m|\leq2^b$ and $1\leq r\leq2^b+2$, then
$\mathrm{Test}(m,r)$ takes $\softO(d^2b)$ bit operations.
\end{lemma}

\begin{proof}
The roots of $q$ are the numbers $((\alpha-m)/r)^{2^N}$, where $\alpha$
runs over the roots of $g$. If the answer is empty,
Lemma~\ref{lem:el-test}(a) gives $|(\alpha-m)/r|^{2^N}>1$ for every
$\alpha$. Otherwise, Lemma~\ref{lem:el-test}(b) gives some $\alpha$ with
$|(\alpha-m)/r|^{2^N}\leq1.5d$, so
$|\alpha-m|\leq r(1.5d)^{1/2^N}\leq2r$.

For the running time, the coefficient of $x^k$ in $p$ is
$\sum_{j\geq k}g_j\binom jkm^{j-k}r^k$, whose absolute value is at most
$2^b\sum_{j=0}^d(|m|+r)^j\leq2^b(d+1)2^{(b+2)d}$. So $p\in\ZZ[x]$ has
coefficients of $L=O(db)$ bits, and \tool{comp:taylor_shift}, applied
over $\ZZ/2^{L+1}\ZZ$, computes $p$ in $\softO(dL)=\softO(d^2b)$ bit
operations. If the coefficients of a polynomial of degree $d$ have
absolute value less than $2^\ell$, then those of $p(x)p(-x)$ have
absolute value less than $(d+1)2^{2\ell}$. Hence the $N$ Graeffe steps
work with integers of at most $2^N(L+\log_2(d+1))=O(db\log d)$ bits.
Each step multiplies two integer polynomials of degree $d$, which takes
$\softO(d^2b)$ bit operations by Kronecker substitution and fast integer
multiplication, and so does the final comparison.
\end{proof}

\paragraph{The search.} We keep a list of integer \emph{centers} at a
\emph{scale}~$h$, where the center $m$ stands for the interval
$[m-h,m+h]$. We start with $h=2^b$ and the single center $0$. At each
scale, we run $\mathrm{Test}(m,h+2)$ for every center $m$ and discard
the centers for which it answers empty. If $h=1$, we output
$\mathcal C'=\{m-1,m,m+1 : m\text{ kept}\}$. Otherwise, we replace each
kept center $m$ by the two centers $m-h/2$ and $m+h/2$, halve $h$, and
repeat.

\begin{proposition}\label{prop:el-search}
The search outputs a set $\mathcal C'\subseteq\ZZ$ that contains every
integer root $r$ of $f$ with $|r|\geq2$. Moreover, $|\mathcal C'|\leq21d$,
\[
 \sum_{c\in\mathcal C'}\bigl\lceil\log_2(2+|c|)\bigr\rceil
 \leq21\bigl(5d+b+\tfrac12\log_2(d+1)\bigr),
\]
and the search takes $\softO(d^3b^2)$ bit operations.
\end{proposition}

\begin{proof}
Let $r$ be an integer root of $f$ with $|r|\geq2$. Every complex root
$z$ of $f$ satisfies $|z|<1+2^b$, since otherwise
$|a_nz^n|\geq|z|^n>2^b\sum_{i<n}|z|^i>\sum_{i<n}|a_iz^i|$, so
$|r|\leq2^b$. By
Lemma~\ref{lem:leading-part} when $d<n$, and trivially when $d=n$, there
is a root $\alpha$ of $g$ with $|r-\alpha|<2$. If $r$ lies in the
interval of a center $m$ at scale $h$, then $|\alpha-m|<h+2$, so by
Lemma~\ref{lem:el-soft} the test does not answer empty and $m$ is kept.
The interval of $m$ is the union of the intervals of the two new
centers $m\mp h/2$ at scale $h/2$. By induction on the scale, $r$ lies
in the interval of a kept center at scale $1$, so $r\in\mathcal C'$.

At scale $h<2^b$ the centers are odd multiples of $h$, so any two of them
are at distance at least $2h$. By Lemma~\ref{lem:el-soft}, every kept
center $m$ has a root $\alpha$ of $g$ with $|\alpha-m|\leq2(h+2)$, and
hence $|\operatorname{Re}\alpha-m|\leq2(h+2)$. For a fixed $\alpha$,
there are at most $4(h+2)/(2h)+1\leq7$ such centers. Since $g$ has at
most $d$ distinct roots, at most $7d$ centers are kept at each scale.
There are $b+1$ scales, so the search runs at most $1+14db$ tests, and
by Lemma~\ref{lem:el-soft} it takes $\softO(d^3b^2)$ bit operations.

At scale $1$ at most $7d$ centers are kept, so $|\mathcal C'|\leq21d$.
Every $c\in\mathcal C'$ is within $1$ of a kept center, which is within
$6$ of a root $\alpha_c$ of $g$. Thus
$2+|c|\leq9+|\alpha_c|\leq10\max(1,|\alpha_c|)$, and
$\lceil\log_2(2+|c|)\rceil\leq5+\log_2\max(1,|\alpha_c|)$. By the
counting above, each root of $g$ serves as $\alpha_c$ for at most $21$
candidates $c$. Summing over the distinct roots $\alpha$ of $g$ and using
Lemma~\ref{lem:landau}, we get
\[
 \sum_{c\in\mathcal C'}\bigl\lceil\log_2(2+|c|)\bigr\rceil
 \leq21\sum_\alpha\bigl(5+\log_2\max(1,|\alpha|)\bigr)
 \leq21\bigl(5d+\log_2M(g)\bigr)
 \leq21\bigl(5d+b+\tfrac12\log_2(d+1)\bigr).\qedhere
\]
\end{proof}

\begin{proof}[Proof of Theorem~\ref{thm:integer-roots}(a)]
Run the search, and evaluate $f$ on
$\mathcal C=\{0,\pm1\}\cup\mathcal C'$ using
Lemma~\ref{lem:multipoint}. By Proposition~\ref{prop:el-search},
$\mathcal C$ contains all integer roots of $f$, $|\mathcal
C|\leq21d+3=O(n)$, and, since $d\leq b$, the bit lengths of its elements
sum to $O(b)$. So Lemma~\ref{lem:multipoint} takes $\softO(nb)$ bit
operations, and the search takes $\softO(d^3b^2)$.
\end{proof}

\section*{Acknowledgments}
\label{sec:acknowledgments}
\addcontentsline{toc}{section}{Acknowledgments}
We thank Prof. Amir Shpilka for helpful advice throughout the work on this paper.

\section*{AI methodology}
\label{sec:ai-methodology}
\addcontentsline{toc}{section}{AI methodology}
We used ChatGPT Work (GPT-6.1-sol, ultra) to help explore the $\widetilde O(nb)$ root-finding algorithm. An initial prompt for a bound linear in $n$ gave an algorithm with time complexity $\widetilde O(nb+b^3)$, which has a cubic dependence on $b$. We next asked it to explore the case where $b$ is larger than $n$, which led the model to suggest Pan's approximate factorization algorithm (Theorem~\ref{thm:pan}), yielding an algorithm with time complexity $\widetilde O(nb+n^2)$. Then we asked it to treat the $n > b$ and $n < b$ regimes separately, which gave the main construction in this paper. Turning this sketch into the present paper required substantial further work. We verified every mathematical detail ourselves, rewrote some of the proofs, and added the related work and algorithmic tools subsections (Sections 1.2 and 3.1). At this stage we also used Claude Opus 5.5 (Anthropic) to improve the wording and to prepare Figure 1. While preparing the final version, we used ChatGPT to analyze Pan's proof of Theorem~\ref{thm:pan}, which led to the concerns described in Remark~\ref{rem:pan}. We used Claude Opus 5.5 to read Pan's paper and confirm these concerns, and to draft Appendix~\ref{sec:elementary}, which we then checked. We verified the correctness and originality of all content, including references. We take full responsibility for the manuscript.

\bibliographystyle{alphaurl}
\bibliography{refs}

\end{document}